\title{On the Density of a Graph and its Blowup}

\author{Asaf Shapira
\thanks{School of Mathematics and College of Computing,
Georgia Institute of Technology, Atlanta GA, 30332. E--mail:
asafico@math.gatech.edu} \and Raphael Yuster\thanks{Department of
Mathematics, University of Haifa, Haifa 31905, Israel. E--mail:
raphy@math.haifa.ac.il} }

\date{}

\documentclass [letterpaper,11pt]{article}
\usepackage{amsfonts}

\newtheorem{theo}{Theorem}
\newtheorem{coro}{Corollary}
\newtheorem{conj}{Conjecture}
\newtheorem{prob}{Problem}

\newtheorem{lemma}{Lemma}[section]
\newtheorem{definition}[lemma]{Definition}
\newtheorem{prop}[lemma]{Proposition}
\newtheorem{comment}[lemma]{Comment}

\newcommand{\qed}{\hspace*{\fill} \rule{7pt}{7pt}}

\newcommand{\ignore}[1]{}

\begin{document}
\maketitle

\begin{abstract}
The theorem of Chung, Graham, and Wilson on quasi-random graphs
asserts that of all graphs with edge density $p$, the random graph
$G(n,p)$ contains the smallest density of copies of
$K_{t,t}$, the complete bipartite graph of size $2t$.
Since $K_{t,t}$ is a $t$-blowup of an edge, the following intriguing open question arises:
Is it true that of all graphs with triangle density $p^3$, the random graph
$G(n,p)$ contains the smallest density of
$K_{t,t,t}$, which is the $t$-blowup of a triangle?

Our main result gives an indication that the answer to the above question is positive
by showing that for {\em some} blowup, the answer must be positive.
More formally we prove that if $G$ has triangle density $p^3$, then there is
some $2 \leq t \leq T(p)$ for which the density of $K_{t,t,t}$ in $G$
is at least $p^{(3+o(1))t^2}$, which (up to the $o(1)$ term) equals the density of $K_{t,t,t}$ in $G(n,p)$. We also
consider the analogous question on skewed blowups, showing that somewhat surprisingly, the behavior there is different.
We also raise several conjectures related to these problems and discuss some applications to other areas.

\end{abstract}

\section{Introduction}\label{intro}

One of the main family of problems studied in extremal graph theory is how does the lack/number of copies
of one graph $H$ in a graph $G$ affects the lack/number of copies of another graph $H'$ in $G$. Perhaps the most
well known problems of this type are Ramsey's Theorem and Tur\'an's Theorem. Our investigation here is concerned
with the relation between the densities of certain fixed graphs in a given graph. Some well known results of this
type are Goodman's Theorem \cite{Goodman} and the Chung-Graham-Wilson Theorem \cite{CGW}. Some recent results of this
type have been obtained recently by Razborov \cite{R} and an abstract investigation of problems of this type was taken
recently by Lov\'asz and Szegedy \cite{LS}. In this paper we introduce an extremal problem of this type, which is
related to some of these well studied problems, and to problems in other areas such as quasi-random graphs and Communication Complexity.

Let us start with some standard notation. Given a graph $H$ on $h$ vertices $v_1,\ldots,v_h$ and a sequence of $h$ positive integers
$a_1,\ldots,a_h$ we denote by $B=H(a_1,\ldots,a_h)$ the $(a_1,\ldots,a_h)$-blowup of $H$ obtained by replacing vertex $v_i$ of $H$ with
an independent set $I_i$ of $a_i$ vertices, and by replacing every edge $(v_i,v_j)$ of $H$ with a complete bipartite graph
connecting the independent sets $I_i$ and $I_j$. For brevity, we will call $B=H(b,\ldots,b)$ the $b$-blowup of $H$, that is,
the blowup in which all vertices are replaced with an independent set of size $b$.
For a fixed graph $H$ and a graph $G$ we denote by $c_H(G)$ the number of copies of $H$
in $G$, or more formally the number of injective mappings from $V(H)$ to $V(G)$ which map edges of $H$ to edges of $G$. For various
reasons, it is usually more convenient to count {\em homomorphisms} from $H$ to $G$ rather than count copies of $H$ in $G$.
Let us then denote this quantity by $Hom_H(G)$, that is, the number of (not necessarily injective)
mappings from $V(H)$ to $V(G)$ which map edges of $H$ to edges of $G$ (allowing two endpoints of an edge to be mapped to the same vertex of $G$). We now let $d_H(G)=Hom_H(G)/n^h$ denote the {\em $H$-density}
of $G$ (or the density of $H$ in $G$). Note that $0 \leq d_H(G) \leq 1$ and we can think of $d_H(G)$ as the probability that a random map $\phi:V(H) \mapsto V(G)$ is a homomorphism. We will also say that a graph on $n$ vertices has edge density $p$ if it has $p{n \choose 2}$ edges. Finally,
let us say that a graph $G$ on $n$ vertices has {\em asymptotic} $H$-density $\gamma$ if $d_H(G)=\gamma \pm o(1)$ where (as usual) the $o(1)$ term
represents a quantity that goes to $0$ when $n$ goes to infinity. For brevity, for the rest of the paper whenever we refer to the $H$-density of a graph we will always refer to the asymptotic $H$-density of $G$.


The main motivation of our investigation here comes from the theory of quasi-random
graphs. The fundamental theorem of quasi-random graphs, the Chung-Graham-Wilson Theorem \cite{CGW} (CGW-Theorem for short),
asserts\footnote{This part of the CGW-Theorems is also implicit in some early works of Erd\H{o}s} that of all graphs with edge density $p$, the random graph
$G(n,p)$ contains the smallest asymptotic density of copies of
$C_4$, the cycle of length $4$. Let $K_{a,b}$ denote the complete bipartite graph on sets of vertices of sizes $a$ and $b$ and note that $K_{a,b}$
is the $(a,b)$-blowup of an edge and that $C_4$ is just $K_{2,2}$. So the CGW-Theorem states
that of all graphs with edge-density $p$, the random graph has the smallest density of the 2-blowup of an edge.
Actually, essentially the same argument as in \cite{CGW} shows that for every $a,b$, of all graphs with edge density
$p$, the random graph contains the smallest density of $K_{a,b}$ (if $a=1$ then any regular graph with edge density $p$
also has this property).

The question we raise in this paper can thus be thought of as an extension of the CGW-Theorem from blowups of an edge, to blowups
of arbitrary graphs. Let us state it explicitly.

\begin{prob}\label{mainprob}
Let $H$ be a fixed graph and set $B=H(a_1,\ldots,a_h)$. Assuming the $d_H(G)=\gamma$, how small can $d_{B}(G)$ be? Furthermore, is it true
that of all graphs with $d_{H}(G)=\gamma$, the (appropriate) random graph $G(n,p)$ has the smallest $d_{B}(G)$?
\end{prob}

As is well-known, the CGW-Theorem further states that if a graph has edge density $p$ and its $C_4$-density is the same
as that of $G(n,p)$ then $G$ must be quasi-random, that is, behave like
$G(n,p)$ in some well defined way (see the excellent survey on quasi-random graphs by
Krivelevich and Sudakov \cite{KSud} for the precise definitions). Again, this result on $C_4$ $~(\equiv K_{2,2})$ can be extended to any $K_{a,b}$ (assuming $K_{a,b}$ is not a star).
We can thus ask the following question that again tries to generalize the result of \cite{CGW} from blowups of an edge to blowups of other graphs $H$:

\begin{prob}\label{mainprobquasi}
Let $H$ be a fixed graph, set $B=H(a_1,\ldots,a_h)$ and suppose $d_H(G)$and $d_{B}(G)$ equal those of $G(n,p)$. Must $G$ be quasi-random?
\end{prob}

As we'll see in this paper, while the result of \cite{CGW} give a positive answer to Problem \ref{mainprob} for all
blowups of an edge, the answer to Problem \ref{mainprob} is negative for some blowups of other graphs, and we conjecture that the answer is
positive for other blowups.
As for Problem \ref{mainprobquasi}, we currently don't have an indication if it has a positive answer for any blowup
of any graph other than blowups of an edge. Hence, we will focus our attention on Problem \ref{mainprob}.
As we will see shortly, Problem \ref{mainprob} seems challenging even for the first non-trivial case of $H$ being
the triangle (denoted $K_3$), so we will mainly consider this special case. To simplify the notation, let us denote by $K_{a,b,c}$ the $(a,b,c)$-blowup of $K_3$. So $K_{2,2,2}$ is the $2$-blowup of the triangle and the question we are interested in is the following:
Suppose the density of triangles in $G$ is $\gamma$. How small can
the density of $B=K_{a,b,c}$ be in $G$? Let us denote by $f_B(\gamma)$ the infimum
of this quantity over all graphs with triangle density $\gamma$. So Problem \ref{mainprob} can be restated as asking
for a bound for the function $f_{B}(\gamma)$.

A simple upper bound for $f_B(\gamma)$ can be obtained by considering the number of triangles
and copies of $K_{a,b,c}$ in the random graph $G(n,\gamma^{1/3})$.
In the other direction, a simple lower bound can be obtained from
the Erd\H{o}s-Simonovits Theorem \cite{ErSim} regarding the number
of copies of complete 3-partite hypergraphs in dense 3-uniform
hypergraphs. These two bounds give the following:

\begin{prop}\label{EasyProp} Let $B=K_{a,b,c}$. Then we have the
following bounds
$$
\gamma^{abc} \leq f_B(\gamma) \leq \gamma^{(ab+bc+ac)/3}~.
$$
\end{prop}

So Problem \ref{mainprob} can be formulated as asking for a characterization of the blowups $B$ for which we have $f_B(\gamma) = \gamma^{(ab+bc+ac)/3}$.
Our main results in this paper, discussed in the next subsections, show that in some cases the lower bound in the above proposition is (essentially) tight, while it seems that for other cases the upper bound gives the correct answer.

We conclude with noting that one can naturally consider the following variant of Problem \ref{mainprob}: Let
$H$ be a fixed graph and let $B'$ be any {\em subgraph} of $H(a_1,\ldots,a_h)$. How small can $f_{B'}(G)$ be if
$f_H(G) = \gamma$? We note that while Problem \ref{mainprob} for the case of $H$ being an edge is well understood (via the CGW-Theorem),
the above variant of Problem \ref{mainprob} is open even when $H$ is an edge. This is the long standing conjecture of Sidorenko \cite{Sid}
and Simonovits \cite{Sim} that states that for every bipartite
graph $B$, the random graph $G(n,p)$ has the smallest $B$-density over all graphs with edge density $p$. We thus focus our attention
on Problem \ref{mainprob}.

\subsection{Balanced blowups and the main results}

When considering the case $B=K_{2,2,2}$, the bounds given by Proposition \ref{EasyProp}
become $\gamma^{8} \leq f_B(\gamma) \leq \gamma^{4}$, and when $B=K_{t,t,t}$ the above bounds
become $\gamma^{t^3} \le f_B(\gamma) \le \gamma^{t^2}$.
The question we are interested in
is, therefore, whether $f_B(\gamma) = \gamma^{t^2}$.
Currently, we cannot obtain any (real) improvement over Proposition \ref{EasyProp}
for the case $B=K_{2,2,2}$, but we are able to show that {\em some} fixed blowup must
have density asymptotically close to $\gamma^{t^2}$, namely, as the density expected
in the random graph. More formally, our main result is the following.

\begin{theo}\label{maintheo}
For every $0 < \gamma, \delta < 1$ there are $N=N(\gamma,\delta)$ and $T=T(\gamma,\delta)$ such
that if $G$ is a graph on $n \geq N$ vertices and its triangle density is $\gamma$, then there is
some $2 \leq t \leq T$ for which the $K_{t,t,t}$-density of $G$
is at least $\gamma^{(1+\delta)t^2}$.
\end{theo}

As we have previously mentioned, we focus our attention on blowups of $K_3$, although
the proof of Theorem \ref{maintheo} extends to blowups of larger complete graphs.
Denote by $K^t_k$ the $t$-blowup of $K_k$. The precise same conclusion of Theorem
\ref{maintheo} holds if $\gamma$ is the density of $K_k$, $T=T(k,\gamma,\delta)$, and $N=N(k,\gamma,\delta)$.
As the proof of Theorem \ref{maintheo} is already quite involved, we only prove it for
triangles, as stated.

Let us state a related result that was recently obtained by Alon \cite{Alon}.

\begin{theo}[Alon \cite{Alon}]\label{theoalon} Set $B=K_{t,t,t}$. Then we have
$$
f_B(\gamma) \geq \gamma^{t^2/\gamma^2}\;.
$$
\end{theo}

Alon's result implies that for any $t \geq 1/\gamma^2$ we can improve upon the lower bound
of Proposition \ref{EasyProp}. Note that Alon's result
does not imply that for large enough $t$, the density of $B=K_{t,t,t}$ gets closer to the density of $B$ in $G(n,p)$.
Alon's argument is based on an idea used by Nikiforov \cite{N} to tackle an Erd\H{o}s-Stone \cite{ErStone} type question.
In Section \ref{open} we observe that a slightly weaker bound can be obtained directly from Nikiforov's result.

Let us finally mention another unexpected motivation for studying $f_B(\gamma)$.
As it turns out, in the case $t=2$ (i.e., when $B=K_{2,2,2}$),
the question of bounding $f_{B}(\gamma)$, was also considered recently (and independently) due to a
different motivation. Barak and Raz observed that improving the
lower bound of $B=K_{2,2,2}$ from $f_B(\gamma) \geq \gamma^{8}$ to
$f_B(\gamma) \geq \gamma^{8-c}$ for some $c>0$ would have certain
non-trivial applications in communication complexity.

\subsection{Skewed blowups}

We now turn our discussion to small skewed-blowups, which somewhat surprisingly, seem to behave quite differently
from the symmetric blowups considered in the previous subsection.
Proposition \ref{EasyProp} implies that when $B=K_{1,1,2}$ we have
$\gamma^{2} \leq f_B(\gamma) \leq \gamma^{5/3}$. In another
independent recent investigation, motivated by an attempt to improve
the bounds in the well-known Triangle Removal Lemma (see Theorem \ref{removallemma}), Trevisan (see \cite{Terry}
page 239) observed that the $\gamma^2$ lower bound for the case
$B=K_{1,1,2}$ can be (slightly) improved:

\begin{theo}\label{theoK112lower} Set $B=K_{1,1,2}$. Then we have the
following bound
$$
f_B(\gamma) \geq \omega(\gamma^{2})\;.
$$
\end{theo}

We note that since the proof of Theorem \ref{theoK112lower} applies the so called triangle removal-lemma (see Theorem
\ref{removallemma}), which, in turn, applies Semer\'edi's Regularity Lemma,
the $\omega(\gamma^{2})$ bound in Theorem \ref{theoK112lower} ``just barely''
beats the simple $\gamma^2$ bound of Proposition \ref{EasyProp}. The bound which the
proof gives is roughly of order $\log^*(1/\gamma)\gamma^2$, and Tao \cite{Terry}
asked if it is possible to improve this bound to something like
$\log\log(1/\gamma)\gamma^2$. While we can not rule out such a
bound, we can still rule out a polynomially better bound by
improving the upper bound of Proposition \ref{EasyProp}.

\begin{theo}\label{theoK112upper} Set $B=K_{1,1,2}$. Then we have the
following bound
$$
f_B(\gamma) \leq \gamma^{2-o(1)}
$$
where the $o(1)$ term goes to $0$ with $\gamma$.
\end{theo}

Observe that Theorem \ref{theoK112upper} implies that the random
graph {\em does not} minimize the density of $K_{1,1,2}$ over all graphs
with a given triangle density. So we see that the answer to Problem \ref{mainprob} is negative
for this case. Also, note that Theorems
\ref{theoK112lower} and \ref{theoK112upper} together determine the
correct exponent of $f_B(\gamma)$ for $B=K_{1,1,2}$. The problem of
determining the correct order of the $o(1)$ terms remains open and
seems challenging.

\begin{comment}
Both Theorems \ref{theoK112lower} and \ref{theoK112upper} were also
obtained independently by N. Alon \cite{Alon}.
\end{comment}

If we consider $B=K_{1,2,2}$, then Proposition \ref{EasyProp} gives
$\gamma^4 \leq f_B(\gamma) \leq \gamma^{8/3}$. Essentially the same
proof as that of Theorem \ref{theoK112lower}, and the same
construction used for the proof of Theorem \ref{theoK112upper}, give
the following improved bounds (we omit the proofs).

\begin{theo}\label{theoK122} Set $B=K_{1,2,2}$. Then we have the following bounds
$$
\omega(\gamma^{4}) \leq f_B(\gamma) \leq \gamma^{3-o(1)}
$$
in which the $o(1)$ term goes to $0$ with $\gamma$, while the
$\omega(1)$ term goes to $\infty$.
\end{theo}

Note that as opposed to the case of $B=K_{1,1,2}$ in which our
bounds determined the correct exponent of $f_B(\gamma)$, in the case
of $B=K_{1,2,2}$ we only know that the correct exponent of
$f_B(\gamma)$ is between $3$ and $4$. Also, $B=K_{1,2,2}$ is another
example of a blowup of $K_3$ for which the answer to Problem \ref{mainprob} is negative.

\subsection{Organization}

The rest of this paper is organized as follows.
In section 2 we focus on large blowups and prove our main result, Theorem \ref{maintheo}.
Our first main tool for the proof of Theorem \ref{maintheo} is the quantitative version of the Erd\H{o}s-Stone theorem, the so called Bollob\'{a}s-Erd\H{o}s-Simonovits theorem \cite{BE,BES}, regarding the size of blowups of $K_r$ in graphs whose density is larger than the Tur\'an density of $K_r$.
Our second main tool is a functional variant of Szemer\'edi's regularity lemma \cite{Sz} due to Alon et al. \cite{AFKS}.
In section 3 we consider small skewed blowups and prove Theorems \ref{theoK112upper} and \ref{theoK112lower}. The proof of these
results apply the so called triangle-removal lemma of Rusza-Szemer\'edi as well as the Rusza-Szemer\'edi graphs.
The final section contains some concluding remarks.

\section{The Density of Large Symmetric Blowups}\label{regularity}

\subsection{Background on the Regularity Lemma}

We start with the basic notions of regularity, some of
the basic applications of regular partitions and state the
regularity lemmas that we use in the proof of Theorem
\ref{maintheo}. See \cite{KS} for a comprehensive survey on the
Regularity Lemma. We start with some basic definitions. For every
pair of nonempty disjoint vertex sets $A$ and $B$ of a graph $G$, we
define $e(A,B)$ to be the number of edges of $G$ between $A$ and
$B$. The {\em edge density} of the pair is defined by
$d(A,B)=e(A,B)/|A||B|$.

\begin{definition}\label{regularpair}\noindent{\bf
($\gamma$-regular pair)} A pair $(A,B)$ is {\em $\gamma$-regular},
if for any two subsets $A' \subseteq A$ and $B' \subseteq B$,
satisfying $|A'| \geq \gamma|A|$ and $|B'| \geq \gamma|B|$, the
inequality $|d(A',B')-d(A,B)| \leq \gamma$ holds.
\end{definition}

Let $G$ be a graph obtained by taking a copy of $K_3$, replacing
every vertex with a sufficiently large independent set, and every
edge with a random bipartite graph. The following well known lemma
shows that if the bipartite graphs are ``sufficiently'' regular, then
$G$ contains the same number of triangles as the random graphs does.
For brevity, let us say that three vertex sets $A,B,C$ are
$\epsilon$-regular if the three pairs $(A,B)$, $(B,C)$ and $(A,C)$
are all $\epsilon$-regular. Several versions of this lemma were
previously proved in papers using the Regularity Lemma. See e.g.
Lemma 4.2 in \cite{F}.

\begin{lemma}\label{cbmsleasy} For every $\zeta$ there is an
$\epsilon=\epsilon_{\ref{cbmsleasy}}(\zeta)$ satisfying the following.
Let $A,B,C$ be pairwise disjoint independent sets of vertices
of size $m$ each that are $\epsilon$-regular and
satisfy $d(A,B)=\alpha_1$, $d(A,C)=\alpha_2$ and
$d(A,C)=\alpha_3$. Then $(A,B,C)$ contain at most
$\left(\alpha_1\alpha_2\alpha_3 + \zeta\right)m^3$ triangles.
\end{lemma}

The following lemma also follows from Lemma 4.2 in \cite{F}.

\begin{lemma}\label{cbmsl} For every $t$ and $\zeta$ there is an
$\epsilon=\epsilon_{\ref{cbmsl}}(t,\zeta)$ such that if $G$ is a $3t$-partite
graph on disjoint vertex sets $A_1,\ldots,A_t$, $B_1,\ldots,B_t$,
$C_1,\ldots,C_t$ of size $m$, and these sets satisfy:
\begin{itemize}
\item $(A_i,B_j,C_k)$ are $\epsilon$-regular for every $1 \leq i,j,k
\leq t$.
\item For every $1 \leq i,j,k \leq t$ we have $d(A_i,B_j) \geq
\alpha_1$, $d(A_i,C_k) \geq \alpha_2$ and $d(B_j,C_k) \geq
\alpha_3$.
\end{itemize}
Then $G$ contains at least $\left(\alpha_1\alpha_2\alpha_3 -
\zeta\right)^{t^2}m^{3t}$ copies of $K_{t,t,t}$ each having precisely one vertex from each partite set.
\end{lemma}

The following lemma is an easy consequence of Lemma \ref{cbmsl}, obtained by taking
$t$ multiple copies of each partite set.
\begin{lemma}\label{oneset} For every $t$ and $\zeta$ there is an
$\epsilon=\epsilon_{\ref{oneset}}(t,\zeta)$ such that if $G$ is a $3$-partite
graph on disjoint vertex sets $A,B,C$ of size $m$
and these sets satisfy:
\begin{itemize}
\item $(A,B,C)$ is $\epsilon$-regular.
\item $d(A,B) \geq \alpha_1$, $d(A,C) \geq \alpha_2$ and $d(B,C) \geq
\alpha_3$.
\end{itemize}
Then $G$ contains at least $\left(\alpha_1\alpha_2\alpha_3 -
\zeta\right)^{t^2}m^{3t}$ distinct homomorphisms of $K_{t,t,t}$.
\end{lemma}

A partition ${\cal A}=\{V_i~|~1\leq i\leq k\}$ of the vertex set of
a graph is called an {\em equipartition} if $|V_i|$ and $|V_{j}|$
differ by no more than $1$ for all $1\leq i < j \leq k$ (so in
particular each $V_i$ has one of two possible sizes). The {\em
order} of an equipartition denotes the number of partition classes
($k$ above). A {\em refinement} of an equipartition ${\cal A}$ is an
equipartition of the form ${\cal B}=\{V_{i,j}~|~1 \leq i \leq k, ~1
\leq j \leq \ell\}$ such that $V_{i,j}$ is a subset of $V_i$ for every
$1 \leq i \leq k$ and $1 \leq j \leq \ell$.

\begin{definition}\label{RegPart}\noindent{\bf ($\gamma$-regular
equipartition)} An equipartition ${\cal B}=\{V_{i}~|~1 \leq i \leq
k\}$ of the vertex set of a graph is called $\gamma$-regular if all
but at most $\gamma k^2$ of the pairs $(V_i,V_{i'})$ are
$\gamma$-regular.
\end{definition}

The Regularity Lemma of Szemer\'edi can be formulated as follows.

\begin{lemma}[\cite{Sz}]\label{SzReg}
For every $m$ and $\gamma>0$ there exists
$T=T_{\ref{SzReg}}(m,\gamma)$ with the following property: If $G$ is
a graph with $n \geq T$ vertices, and ${\cal A}$ is an equipartition
of the vertex set of $G$ of order at most $m$, then there exists a
refinement ${\cal B}$ of ${\cal A}$ of order $k$, where $m \leq k
\leq T$ and ${\cal B}$ is $\gamma$-regular.
\end{lemma}

Our main tool in the proof of Theorem \ref{maintheo} is Lemma
\ref{NewReg1} below, proved in \cite{AFKS}. This lemma can be
considered a strengthening of Lemma \ref{SzReg}, as it guarantees
the existence of an equipartition and a refinement of this
equipartition that poses stronger properties compared to those of
the standard $\gamma$-regular equipartition. This stronger notion is
defined below.

\begin{definition}\label{ERegPart}\noindent{\bf (${\cal E}$-regular
equipartition)} For a function ${\cal E}(r): \mathbb{N} \mapsto
(0,1)$, a pair of equipartitions ${\cal A}=\{V_{i}~|~1 \leq i \leq
k\}$ and its refinement ${\cal B}=\{V_{i,j}~|~1 \leq i \leq k, ~1
\leq j \leq \ell \}$, where $V_{i,j} \subset V_i$ for all $i,j$, are
said to be ${\cal E}$-regular if
\begin{enumerate}

\item All but at most ${\cal E}(0)k^2$ of the pairs $(V_i,V_j)$ are ${\cal
E}(0)$-regular.

\item For all $1 \leq i < i' \leq k$, for all $1 \leq j,j' \leq \ell$
but at most ${\cal E}(k)\ell^2$ of them, the pair $(V_{i,j},V_{i',j'})$
is ${\cal E}(k)$-regular.

\item All $1 \leq i < i' \leq k$ but at most ${\cal
E}(0)k^2$ of them are such that for all $1 \leq j,j' \leq \ell$ but at
most ${\cal E}(0)\ell^2$ of them $|d(V_i,V_{i'})-d(V_{i,j},V_{i',j'})|
< {\cal E}(0)$ holds.
\end{enumerate}
\end{definition}

It will be very important for what follows to observe that in
Definition \ref{ERegPart} we may use an arbitrary {\em function}
rather than a fixed $\gamma$ as in Definition \ref{RegPart} (such
functions will be denoted by ${\cal E}$ throughout the paper). The
following is one of the main results of \cite{AFKS}.

\begin{lemma}\noindent{\bf(\cite{AFKS})}\label{NewReg1}
For any integer $m$ and function ${\cal E}(r): \mathbb{N} \mapsto
(0,1)$ there is $S=S_{\ref{NewReg1}}(m,{\cal E})$ such that any
graph on at least $S$ vertices has an ${\cal E}$-regular
equipartition ${\cal A}$, ${\cal B}$ where $|{\cal A}|=k \geq m$ and
$|{\cal B}|=k\ell \leq S$.
\end{lemma}

\subsection{Main Idea and Main Obstacle}

Let us describe the main intuition behind the proof of Theorem \ref{maintheo},
and where its naive implementation fails. Recall that Lemma \ref{cbmsleasy} says
that an $\epsilon$-regular triple contains the ``correct'' number of
triangles we expect to find in a ``truly'' random graph with the same density. So given a graph with triangle
density $\gamma$, we can apply the Regularity Lemma with (say)
$\epsilon=\gamma$. Suppose we get a partition into $k$ sets, for
(say) some $k \leq T_{\ref{SzReg}}(\gamma, 1/\gamma^2)$. So the situation now is that
the number of triangles spanned by any triple $V_i,V_j,V_k$ is more
or less determined by the densities between the sets. Since $G$ has
triangle density $\gamma$, we get (by averaging) that there must be
some triple $V_i,V_j,V_k$ whose triangle density is also close to
being at least $\gamma$. Suppose the densities between $V_i,V_j,V_k$ are
$\alpha_1$, $\alpha_2$ and $\alpha_3$. Since the number of triangles
between $V_i,V_j,V_k$ is determined by the densities connecting them,
we get that $\alpha_1\alpha_2\alpha_3$ is close to $\gamma$. Now, if
$\epsilon$ is small enough, then we can also apply Lemma \ref{cbmsl}
on $V_i,V_j,V_k$ in order to infer that they contain close to
$(n/k)^{3t}\alpha_1^{t^2}\alpha_2^{t^2}\alpha_3^{t^2}$ copies of
$K_{t,t,t}$. Hence, by the above consideration, this
number is close to $(n/k)^{3t}\gamma^{t^2}$. Now, since for
large enough $t \geq t(k)$ we have $(n/k)^{3t}\gamma^{t^2} =
n^{3t}\gamma^{(1+o(1))t^2}$ we can choose a large enough $t=t(k)$ to
get the desired result. Since $k$ is bounded by a function of
$\gamma$ so is $t$.

The reason why the above argument fails is that in order to apply
Lemma \ref{cbmsl} with a given $t$, the value of $\epsilon$ in the
$\epsilon$-regular partition needs to depend on $t$. So we arrive at
a circular situation in which $\epsilon$ needs to be small enough in
terms of $t$ (to allow us to apply Lemma \ref{cbmsl}), and $t$ needs
to be large enough in terms of $\epsilon$ (to allow us to infer that
$(n/k)^{3t}\gamma^{t^2} = n^{3t}\gamma^{(1+o(1))t^2}$).

We overcome the above problem by applying Lemma \ref{NewReg1} which
more or less allows us to find a partition which is $f(k)$-regular
where $k$ is the number of partition classes. However, this is an over
simplification of this result (as can be seen from Definition \ref{ERegPart}), and our proof
requires several other ingredients that enable us to apply Lemma \ref{NewReg1}.
Most notably, we need to use a classic result of Bollob\'{a}s, Erd\H{o}s and Simonovits \cite{BE,BES} and adjust it
to our setting.

\subsection{Some preliminary lemmas}

We now turn to discuss two simple (yet crucial) lemmas that will be later used in the proof
of Theorem \ref{maintheo}. Let us recall that
Tur\'an's theorem asserts that every graph with edge density larger
than $1-\frac{1}{r-1}$ contains a copy of $K_r$, the complete graph
on $r$ vertices. The Erd\H{o}s-Stone theorem strengthens this result
by asserting that if the edge density is larger than
$1-\frac{1}{r-1}$, then the graph actually contains a blowup of
$K_r$. More precisely, there is a function $f(n,c,r)$ such that
every $n$-vertex graph with edge density $1-\frac{1}{r-1}+\beta$
contains an $f(n,\beta,r)$-blowup of $K_r$. The determination of the
growth rate of $f(n,\beta,r)$ received a lot of attention until Bollob\'{a}s,
Erd\H{o}s and Simonovits \cite{BE,BES} determined that for fixed $\beta$ and $r$ we have
$f(n,\beta,r)=\Theta(\log n)$. See \cite{N} for a short proof of
this result and for related results and references. As it turns out, the
bound $\Theta(\log n)$ will be crucial for our proof (a bound like $\log^{1-\epsilon} n$ would not be useful for us). Let us state an
equivalent formulation of this result for the particular choice of
$r=3$ and $\beta=1/24$.

\begin{theo}[Bollob\'{a}s-Erd\H{o}s-Simonovits \cite{BE,BES}]\label{BESS} There is
an absolute constant $c$, such that every graph on at least $c^{t}$
vertices and edge density at least $13/24$ contains a copy of $K_{t,t,t}$.
\end{theo}
As a $3$-partite graph with edge density at least $7/8$ between any two parts has
overall density greater than $13/24$ we have:
\begin{coro}
\label{coro-ttt}
There is an absolute constant $C$, such that every 3-partite graph with parts
of equal size $C^{t}$ and edge density at least $7/8$ between any two parts, contains a copy of
$K_{t,t,t}$.
\end{coro}


We will need the following lemma guaranteeing many copies of a large blowup of $K_3$.

\begin{lemma}\label{pack}
If $G$ is a 3-partite graph on vertex sets $X$, $Y$ and $Z$ of equal size $m$, and the
three densities $d(X,Y)$, $d(X,Z)$ and $d(Y,Z)$ are all at least
$15/16$, then $G$ contains at least $\lfloor m^{3t}/C^{3t^2}\rfloor$ copies of $K_{t,t,t}$.
\end{lemma}

\paragraph{Proof:}
Let $C$ be the constant of Corollary \ref{coro-ttt}.
If $m < C^t$ there is nothing to prove (as $\lfloor m^{3t}/C^{3t^2} \rfloor =0$) so let us
assume that $m \geq C^t$.
We first claim that at least $1/2$
of the graphs spanned by three sets of vertices $X' \subseteq X$,
$Y' \subseteq Y$, $Z' \subseteq Z$, where
$|X'|=|Y'|=|Z'|=C^t$, have edge density at least $7/8$.
Indeed, suppose we randomly pick the sets $X'$,
$Y'$ and $Z'$. The expected density of non-edges between $(X',Y')$,
$(X',Z')$ and $(Y',Z')$ is $1/16$, so by Markov's
inequality, with probability at least $1/2$ this density is at most
$1/8$.

By Corollary \ref{coro-ttt}, every graph of size at least $C^t$,
whose edge density is at least $7/8$, contains a copy of
$K_{t,t,t}$. So by the above consideration, at least half of the
${m \choose {C^t}}^3$ choices of $A',B',C'$ contain a $K_{t,t,t}$.
Since each such $K_{t,t,t}$ is counted ${{m-t} \choose {C^t-t}}^3$
times, we have that the number of distinct copies of $K_{t,t,t}$ in $G$ is at
least
$$
\frac12{m \choose {C^t}}^3 / {{m-t} \choose {C^t-t}}^3
\geq m^{3t}/C^{3t^2}\;.
$$

$\qed$

\bigskip

The proof of Theorem \ref{maintheo} we give in the next subsection only covers the case of $\gamma  \ll \delta$. As the following lemma shows,
we can then ``lift'' this result to arbitrary $0< \gamma,\delta <1$.

\begin{lemma}\label{landau} If Theorem \ref{maintheo} holds for every $\delta >0$ and every
small enough $\gamma < \gamma_0(\delta)$, then it also holds for every $0 < \gamma,\delta < 1$.
\end{lemma}

\paragraph{Proof:} Assume to the contrary that there exist a $\delta >0$, a $\gamma \geq \gamma_0(\delta)$ and
arbitrarily large graphs with triangle-density $\gamma$ in which the $K_{t,t,t}$-density is smaller than $\gamma^{(1+\delta)t^2}$ for every
$2 \leq t \leq T(\gamma^2_0(\delta),\delta)$. Let $G$ be one such graph on at least $N(\gamma^2_0(\delta),\delta)$ vertices.
For an integer $k$ let $G^{\otimes k}$ be the $k^{th}$ tensor product of $G$,
that is, the graph whose vertices are sequences of $k$ (not necessarily distinct) vertices of $G$, and where vertex
$v=(v_1,\ldots,v_k)$ is connected to vertex $u=(u_1,\ldots,u_k)$ if and only if $v_i$ is connected to $u_i$ for every $1 \leq i \leq k$.
The key observation is that for every graph $H$, if the $H$-density of $G$ is $\gamma$ then the $H$-density of $G^{\otimes k}$ is $\gamma^k$.
Let then $k$ be the smallest integer satisfying $\gamma^k < \gamma_0(\delta)$ and note that in this case we have
$\gamma^2_0(\delta) \leq  \gamma^k < \gamma_0(\delta)$. We thus get that $G^{\otimes k}$ is a graph on at least
$N(\gamma^2_0(\delta),\delta) \geq N(\gamma^k,\delta)$ vertices, with triangle density $\gamma^k$ and for all $2 \leq t \leq T(\gamma^k,\delta) \leq T(\gamma^2_0(\delta),\delta)$ its $K_{t,t,t}$ density is smaller than $\gamma^{k(1+\delta)t^2}$, which contradicts the assumption of the lemma. $\qed$

\subsection{Proof of Theorem \ref{maintheo}}\label{subsecmain}

We prove the theorem for every $0 < \delta < 1$ and for every $0 < \gamma < 1$ which is small enough so that
\begin{equation}\label{ChooseGamma}
\gamma < \left(\frac{1}{128C^3}\right)^{2/\delta}\;,
\end{equation}
where $C$ is the absolute constant from Lemma \ref{pack}.
By Lemma \ref{landau} this will establish the theorem for all $0 < \delta, \gamma < 1$.

For a given positive integer $r$, let $t=t(r,\delta,\gamma)$ be a
large enough integer such that
\begin{equation}\label{keyeq}
\frac{1}{r^{3t}}\left(\frac{\gamma}{64}\right)^{t^2} \geq
2C^{3t^2}\gamma^{(1+\delta)t^2}
\end{equation}
holds. Since we assume that $\gamma$ and $\delta$ satisfy (\ref{ChooseGamma}), it is enough to make sure that $t$ satisfies
$$
\frac{1}{r^{3t}} \geq \gamma^{\frac12\delta t^2}\;,
$$
hence we can take
\begin{equation}\label{definet}
t(r,\delta,\gamma)=\max\{2,~\frac{6\log r}{\delta \log \frac{1}{\gamma}}\}\;.
\end{equation}
We now define a function ${\cal E}(r)$ as follows:
\begin{equation}\label{EqBig1}
{\cal E}(r)=
\left\{%
\begin{array}{ll}
    \min\{\frac{1}{16}, ~\gamma/30,~\epsilon_{\ref{cbmsleasy}}(\gamma/4)\}, & \hbox{$r=0$}
    \\~\\
    \min\{\frac{1}{16}, ~\epsilon_{\ref{oneset}}(t(r,\delta,\gamma),\gamma/64)    ,~\epsilon_{\ref{cbmsl}}(t(r,\delta,\gamma),\gamma/64) \}, & \hbox{$r \geq 1$}\;. \\
\end{array}%
\right.
\end{equation}

Given $\gamma$ and $\delta$ let ${\cal E}(r)$ be the function defined above.
Set $m=30/\gamma$ and let $S=S_{\ref{NewReg1}}(m,{\cal E})$ be the constant from
Lemma \ref{NewReg1}. Given a graph $G$ on $n \geq S$ vertices and parameters
$\gamma$ and $\delta$, we apply Lemma \ref{NewReg1} with $m=30/\gamma$ and with
the function ${\cal E}(r)$ defined above. The lemma returns an
${\cal E}$-regular partition consisting of an equipartition ${\cal
A}=\{V_{i}~|~1 \leq i \leq k\}$ and a refinement ${\cal
B}=\{V_{i,j}~|~1 \leq i \leq k, ~1 \leq j \leq \ell \}$, where $k\ell \leq
S(m,{\cal E})$ and $k \ge m$. Note that $S$ depends only on $\delta$ and
$\gamma$.

We now remove from $G$ any edge whose endpoints belong to
the same set $V_i$. We thus remove at most $n^2/(2m) < \frac{\gamma}{60}n^2$
edges. We also remove any edge connecting pairs $(V_i,V_j)$ that are
not ${\cal E}(0)$-regular. The first property of an ${\cal
E}$-regular partition guarantees that we thus remove at most ${\cal
E}(0)n^2 \leq \frac{\gamma}{30}n^2$ edges. We also remove any edge
connecting a pair $(V_i,V_j)$ for which there are more that ${\cal
E}(0)\ell^2$ pairs $i',j'$ which do not satisfy
$|d(V_i,V_{j})-d(V_{i,i'},V_{j,j'})| < {\cal E}(0)$. By the third
property of an ${\cal E}$-partition we infer that we thus remove
at most $\frac{\gamma}{30}n^2$ edges. All together we have removed
less than $\frac{\gamma}{12} n^2$ edges and so we have destroyed at most
$\frac{\gamma}{2} n^3$ triangles in $G$ (recall that we are counting homomorphisms so each triangle is counted $6$ times).
And so the new graph we obtain has triangle density at least $\gamma/2$. Let us call this new graph
$G'$.

As $G'$ has triangle density at least $\gamma/2$, we get (by
averaging) that there must be three sets $(V_i,V_j,V_k)$ that
contain at least $\frac12\gamma(n/k)^3$ triangles with one vertex in
each of the sets $V_i,V_j,V_k$ (we are using $k$ as both an index and as the number of parts in the partition, but there is no confusion). For what follows, let us set
$\alpha_1=d(V_i,V_j)$, $\alpha_2=d(V_i,V_k)$ and
$\alpha_3=d(V_j,V_k)$. Because we have removed edges between
non-${\cal E}(0)$-regular pairs, we get that $(V_i,V_j,V_k)$ must be
${\cal E}(0)$-regular. Letting $\Delta$ denote the number of
triangles spanned by $(V_i,V_j,V_k)$ we see that as ${\cal E}(0)
\leq \epsilon_{\ref{cbmsleasy}}(\gamma/4)$, we can apply Lemma
\ref{cbmsleasy} on $(V_i,V_j,V_k)$ to conclude that
$$
\frac12\gamma\left(\frac{n}{k}\right)^3\leq \Delta \leq
(\alpha_1\alpha_2\alpha_3+\frac14\gamma)\left(\frac{n}{k}\right)^3\;,
$$
implying that

\begin{equation}\label{boundgamma}
\alpha_1\alpha_2\alpha_3 \geq \frac14\gamma\;.
\end{equation}

Let us say that a $3s$-tuple (where $s$ is any positive integer) $1 \leq i_1 <\cdots<i_s \leq \ell$, $1
\leq j_1 < \cdots < j_s \leq \ell$, $1 \leq k_1 \cdots < k_s \leq \ell$ is
{\em good} if it satisfies the following four properties:
\begin{enumerate}
\item For every $i_a,j_b,k_c$ we have that
$(V_{i,i_a},V_{j,j_b},V_{k,k_c})$ are ${\cal E}(k)$-regular.
\item For every $i_a,j_b$ we have
$d(V_{i,i_a},V_{j,j_b}) \geq \alpha_1 - {\cal E}(0) \geq
\alpha_1-\frac18\gamma \geq \frac12\alpha_1$.
\item For every $i_a,k_c$ we have
$d(V_{i,i_a},V_{k,k_c}) \geq \alpha_2 - {\cal E}(0) \geq
\alpha_2 - \frac18\gamma \geq \frac12\alpha_2$.
\item For every $j_b,k_c$ we have
$d(V_{j,j_b},V_{k,k_c}) \geq \alpha_3 - {\cal E}(0) \geq
\alpha_3 - \frac18\gamma \geq \frac12\alpha_3$.
\end{enumerate}

Suppose $i_1,\ldots,i_t$, $j_1,\ldots,j_t$, $k_1,\ldots,k_t$ is a
good $3t$-tuple. Then the definition of ${\cal E}$ (via the function
$\epsilon_{\ref{cbmsl}}(t,\zeta)$ from Lemma \ref{cbmsl}) and the first property
of a good $3t$-tuple, guarantee that we can apply Lemma \ref{cbmsl}
on $V_{i,i_1},\ldots,V_{i,i_t}$, $V_{j,j_1},\ldots,V_{j,j_t}$,
$V_{k,k_1},\ldots,V_{k,k_t}$, to conclude that they have at least
$$
\left(\frac{n}{kl}\right)^{3t}\left(\frac18\alpha_1\alpha_2\alpha_3-
\frac{1}{64}\gamma\right)^{t^2}
\geq
\left(\frac{n}{kl}\right)^{3t}\left(\frac{\gamma}{64}\right)^{t^2}
$$
copies of $K_{t,t,t}$, where we have also used (\ref{boundgamma}).
Our choice of $t=t(k,\delta,\gamma)$ in (\ref{definet}) guarantees
(via (\ref{keyeq})) that the number of copies of $K_{t,t,t}$ in a
good $3t$-tuple is at least
\begin{equation}\label{firstLB}
\left(\frac{n}{k\ell}\right)^{3t}\left(\frac{\gamma}{64}\right)^{t^2}
\geq 2C^{3t^2}\left(\frac{n}{\ell}\right)^{3t}\gamma^{(1+\delta)t^2}\;.
\end{equation}

But how can we be certain that a good $3t$-tuple exists? And if they do, how many are there?
We first consider the case $\ell \ge C^t$.
Let us now recall that ${\cal E}(r) \leq \frac{1}{16}$ for every $r
\geq 0$ and so the second and third properties of a ${\cal
E}$-regular partition guarantee that at least $\frac{15}{16}\ell^2$ of the
choices $1 \leq i',j' \leq \ell$ are such that $(V_{i,i'},V_{j,j'})$ is
${\cal E}(k)$-regular and satisfies $|d(V_i,V_j) -
d(V_{i,i'},V_{j,j'})| \leq {\cal E}(0)$. The same holds with respect
to the other two pairs $(V_{j},V_{k})$ and $(V_{i},V_{k})$.
Therefore, by Lemma \ref{pack}, the sets $V_i,V_j,V_k$ contain at least
$\lfloor \ell^{3t}/C^{3t^2} \rfloor \ge 0.5\ell^{3t}/C^{3t^2}$ choices of good $3t$-tuples. Hence, combining this
with (\ref{firstLB}) we infer that the number of copies of
$K_{t,t,t}$ spanned by $(V_i,V_j,V_k)$ is at least
$$
\frac{\ell^{3t}}{2C^{3t^2}} \cdot
2C^{3t^2}\left(\frac{n}{\ell}\right)^{3t}\gamma^{(1+\delta)t^2}=
n^{3t}\gamma^{(1+\delta)t^2}\;,
$$
implying that the density of $K_{t,t,t}$ in $G'$ (and so also in
$G$) is at least $\gamma^{(1+\delta)t^2}$.

We now consider the case $\ell < C^t$.
Assume that in this case we can find just {\em one} good $3$-tuple.
Then definition of ${\cal E}$ (via the function
$\epsilon_{\ref{oneset}}(t,\zeta)$ from Lemma \ref{oneset}) and the first property
of a good $3$-tuple, guarantee that we can apply
Lemma \ref{oneset} on this $3$-tuple, to conclude that it has at least
$$
\left(\frac{n}{kl}\right)^{3t}\left(\frac18\alpha_1\alpha_2\alpha_3-
\frac{\gamma}{64}\right)^{t^2}
\geq
\left(\frac{n}{kl}\right)^{3t}\left(\frac{\gamma}{64}\right)^{t^2}
$$
distinct homomorphisms of $K_{t,t,t}$.
Our choice of $t=t(k,\delta,\gamma)$ in (\ref{definet}) guarantees
(via (\ref{keyeq})) that the number of homomorphisms of $K_{t,t,t}$ in a
good $3$-tuple is at least
$$
\left(\frac{n}{k\ell}\right)^{3t}\left(\frac{\gamma}{64}\right)^{t^2}
\geq 2C^{3t^2}\left(\frac{n}{\ell}\right)^{3t}\gamma^{(1+\delta)t^2} \geq
n^{3t}\gamma^{(1+\delta)t^2}\;,
$$
implying that the density of $K_{t,t,t}$ in $G'$ (and so also in
$G$) is at least $\gamma^{(1+\delta)t^2}$. To see that a single good
$3$-tuple $i_1$, $j_1$, $k_1$ exists, consider picking $i_1$, $j_1$ and $k_1$
randomly and uniformly from $[\ell]$. Since ${\cal E}(k),{\cal E}(0) \leq \frac{1}{16}$
we infer that with positive probability $i_1$, $j_1$ and $k_1$ will satisfy the four
properties of a good $3$-tuple, so a good $3$-tuple exists.

Finally, note that since $k \leq S$ we see that $k$ is upper
bounded by some function of $\gamma$ and $\delta$. As
$t=t(k,\delta,\gamma)$ is chosen in (\ref{definet}) we see that $2 \leq t \leq T(\gamma,\delta)$
for some function $T(\gamma,\delta)$ and so the proof is complete. $\qed$

\section{The Density of Small Skewed Blowups}

In this section we focus our attention on small skewed blowups of $K_3$.
We start with that proof of Theorem \ref{theoK112upper} in which we will apply
the following well known result of Ruzsa and Szemer\'edi \cite{RuS}. For
completeness, we include the short proof.

\begin{theo}[\cite{RuS}]\label{RuzSzeBeh}
Suppose $S \subseteq [n]$ is a set of integers containing no 3-term
arithmetic progression. Then there is a graph $G=(V,E)$ with
$|V|=6n$ and $|E|=3n|S|$, whose edges can be (uniquely) partitioned
into $n|S|$ edge disjoint triangles. Furthermore, $G$ contains no
other triangles.
\end{theo}

\paragraph{Proof:} We define a 3-partite
graph $G$ on vertex sets $A$, $B$ and $C$, of sizes $n$, $2n$ and
$3n$ respectively, where we think of the vertices of $A$, $B$ and
$C$ as representing the sets of integers $[n]$, $[2n]$ and $[3n]$. For every
$1 \leq i \leq n$ and $s \in S$ we put a triangle $T_{i,s}$ in $G$
containing the vertices $i \in A$, $i+s \in B$ and $i+2s \in C$. It
is easy to see that the above $n|S|$ triangles are edge disjoint,
because every edge determines $i$ and $s$. To see that $G$ does not
contain any more triangles, let us observe that $G$ can only contain
a triangle with one vertex in each set. If the vertices of this
triangle are $a \in A$, $b \in B$ and $c \in C$, then we must have
$b=a+s_1$ for some $s_1 \in S$, $c=b+s_2=a+s_1+s_2$ for some $s_2
\in S$, and $a=c-2s_3=a+s_1+s_2-2s_3$ for some $s_3$. This means
that $s_1,s_2,s_3 \in S$ form an arithmetic progression, but because
$S$ is free of 3-term arithmetic progressions it must be the case
that $s_1=s_2=s_3$ implying that this triangle is one of the
triangles $T_{i,s}$ defined above. $\qed$

\bigskip

For the proof of Theorem \ref{theoK112upper} we will need to combine
Theorem \ref{RuzSzeBeh} with the following well known result of
Behrend \cite{B}.

\begin{theo}[Behrend \cite{B}]\label{Behrend}
For every $n$, there exists $S \subseteq [n]$ of size
$n/8^{\sqrt{\log n}}=n^{1-o(1)}$ containing no 3-term arithmetic
progression.
\end{theo}

\paragraph{Proof of Theorem \ref{theoK112upper}:}
Let $m$ be any integer and let $S$ be a $3AP$-free subset of $[m]$
of size $m/8^{\sqrt{\log m}}$ as guaranteed by Theorem
\ref{Behrend}. Let $G'$ be the graph of Theorem \ref{RuzSzeBeh} when
using $[m]$ and the above set $S$. Finally, let $G$ be an $n/6m$
blowup of $G'$, that is, the graph obtained by replacing every
vertex $v$ of $G'$ with an independent set $I_v$ of size $n/6m$, and
replacing every edge $(u,v)$ of $G'$ with a complete bipartite graph
connecting $I_u$ and $I_v$. Observe that $G$ has $n$ vertices, and
that each triangle in $G'$ gives rise to $(n/6m)^3$ triangles in
$G$. Hence, the number of ways to map a triangle into $G$ is
$$
6m|S|\left(\frac{n}{6m}\right)^3 = \frac{n^3}{6^{2}m8^{\sqrt{\log
m}}}\;
$$
(recall that there are six ways to map a labeled triangle into a
triangle of $G$). The crucial observation is that because all the
triangles in $G'$ are edge disjoint, the only copies of $K_{1,1,2}$
in $G$ are those that are formed by picking 4 vertices from the sets
$I_a$, $I_b$ and $I_c$ for which $a$, $b$ and $c$ formed a triangle
in $G'$. This means that the number of ways to map a $K_{1,1,2}$
into $G$ is
$$
12m|S|\left(\frac{n}{6m}\right)^2{{\frac{n}{6m}} \choose 2} \le
\frac{n^4}{6^{3}m^28^{\sqrt{\log m}}}\;.
$$
Now setting
$$
\gamma=\frac{1}{6^{2}m8^{\sqrt{\log m}}}
$$
we see that the density of triangles in $G$ is $\gamma$, while the
density of $K_{1,1,2}$ in $G$ is at most
$$
\frac{1}{6^{3}m^28^{\sqrt{\log m}}}=\gamma^{2}2^{c\sqrt{\log
1/\gamma}}=\gamma^{2-o(1)}\;,
$$
for some absolute constant $c$, thus completing the proof. $\qed$

\bigskip

For completeness, we reproduce the short proof of Theorem \ref{theoK112lower}.
We will need the so called ``triangle removal lemma'' of \cite{RuS}:

\begin{theo}[\cite{RuS}]\label{removallemma}
If $G$ is an $n$ vertex graph from which one should remove at least
$\epsilon n^2$ edges in order to destroy all triangles, then $G$
contains at least $f(\epsilon) n^3$ triangles.
\end{theo}

\paragraph{Proof of Theorem \ref{theoK112lower}:}
Suppose $G$ has $\gamma n^3$ triangles. Then by Theorem
\ref{removallemma} we know that $G$ contains a set of edges $F$ of
size at most $f(\gamma)n^2$, the removal of which makes $G$
triangle-free, where $f(\gamma)=o(1)$. For each edge $e \in E(G)$
let $t(e)$ be the number of triangles in $G$ containing $e$ as one
of their edges. Observe that a copy of $K_{1,1,2}$ is obtained by
taking two triangles sharing an edge. Also, as the removal of edges
in $F$ makes $G$ triangle-free, every triangle in $G$ has an edge of
$F$ as one if its edges. Therefore, by Cauchy-Schwartz we have that
the number of copies of $K_{1,1,2}$ in $G$ is
$$
\sum_{e \in F}t(e)^2 \geq \frac{1}{|F|}\left(\sum_{e \in F}t(e)
\right)^2 \geq \frac{1}{|F|}\gamma^2n^6 \geq
\frac{1}{f(\gamma)}\gamma^2n^4\;,
$$
implying the desired result with $1/f(\gamma)$ being the $\omega(1)$
term in the statement of the theorem. $\qed$

\section{Concluding Remarks and Open Problems}\label{open}

\begin{itemize}

\item Our main result given in Theorem \ref{maintheo} states that in any graph $G$ with $K_3$-density
$\gamma$, there is {\em some} $t$ for which the $K_{t,t,t}$-density in $G$ is almost as large as the
$K_{t,t,t}$-density in a random graph with the same triangle density. This motivates us to raise the following conjecture, stating that the
upper bound in Proposition \ref{EasyProp} gives the correct order of $f_{B}(\gamma)$ and hence that the answer to Problem \ref{mainprob} is positive for balanced blow-ups of $K_3$.

\begin{conj}\label{ConjBalance}
Let $t \geq 2$ and set $B=K_{t,t,t}$. Then
$$
f_{B}(\gamma)=\gamma^{t^2}\;.
$$
\end{conj}

We remind the reader of our remark in Section \ref{intro}, that any (polynomial) improvement over the lower bound
of Proposition \ref{EasyProp} would have interesting applications in theoretical computer science.

\item Theorems \ref{theoK112upper} and the upper bound of \ref{theoK122} show that when considering
the skewed blowups $B=K_{1,1,2}$ or $B=K_{1,2,2}$ the random graph does not minimize the density of $B$.
This motivates us to raise the following conjecture.

\begin{conj}\label{ConjSkew}
If $a,b,c$ are not all equal, then there is $c>0$ such that
$$
f_{B}(\gamma) \leq \gamma^{\frac13(ab+bc+ac)+c}\;.
$$
\end{conj}

The construction we used in order to prove Theorems \ref{theoK112upper} and \ref{theoK122} can be used
to verify Conjecture \ref{ConjSkew} for other skewed blowups. For example, for every $B=K_{1,1,t}$ it
establishes that $f_{B} \leq \gamma^{t-o(1)}$ which matches the lower bound in Proposition \ref{EasyProp} up to the $o(1)$ term.
It may be possible to modify this construction in order to establish Conjecture \ref{ConjSkew}.

\item As we have mentioned in Section \ref{intro}, Alon \cite{Alon} has recently shown that
if the triangle density of $G$ is $\gamma$ then its $K_{t,t,t}$-density is at least $\gamma^{O(t^2/\gamma^2)}$.
We now show that a slightly weaker bound can be derived directly from a recent result of Nikiforov \cite{N}.

\begin{theo}\label{NT} If a graph has triangle density $\gamma$, then its
$K_{t,t,t}$-density is at least $2^{-O(t^2/\gamma^3)}$.
\end{theo}

\paragraph{Proof (sketch):} By a result of Nikiforov \cite{N},
a graph with triangle density $\gamma$ has a
$K_{t,t,t}$ with $t=\gamma^3\log n$. Or in other words, every graph on at
least $2^{t/\gamma^3}$ vertices, whose triangle density is $\gamma$,
has a copy of $K_{t,t,t}$. As in the proof of Lemma \ref{pack}, if a
graph has triangle density $\gamma$, then most subsets of vertices
of size $2^{t/\gamma^3}$ have (roughly) the same density, so they
contain a $K_{t,t,t}$. We thus get that $G$ has $\frac12{n \choose
2^{t/\gamma^3}}$ sets which contain a $K_{t,t,t}$ and since each
$K_{t,t,t}$ is counted ${n-3t \choose 2^{t/\gamma^3}-3t}$ times we
get that $G$ has $n^{3t}/2^{O(t^2/\gamma^3)}$ distinct copies of
$K_{t,t,t}$. $\qed$

\bigskip

We note that although this is not stated explicitly in \cite{N}, Nikiforov's arguments actually shows
that a graph of triangle density $\gamma$ has a $K_{t,t,t}$ with $t=\gamma^2\log_{1/\gamma} n$ and so the argument above can actually
give Alon's result.

\item Observe that in a random graph $G(n,\gamma^{1/3})$, whose triangle
density is $\gamma$, we expect to find a $K_{t,t,t}$ with
$t=c\log_{1/\gamma} n$ for some absolute constant $c$. It seems very interesting to try and
improve Nikiforov's result \cite{N} mentioned above by showing the following:

\begin{conj}
There is an absolute constant $c>0$, such that if a graph $G$ has triangle density $\gamma$, then $G$ has a $K_{t,t,t}$
of size $t=c\log_{1/\gamma} n$.
\end{conj}

Besides being an interesting problem on its own, we note that such an
improved bound, together with the argument we gave in the proof of
Theorem \ref{NT}, would imply that if the triangle density of a
graph is $\gamma$, then its $K_{t,t,t}$ density is at least
$\gamma^{O(t^2)}$, which would come close to establishing Conjecture \ref{ConjBalance}.

\end{itemize}

\paragraph{Acknowledgement:} We would like to thank Noga Alon, Guy Kindler and Benny Sudakov
for helpful discussions related to this paper.

\end{document}